\documentclass[12pt]{amsart}
 
\usepackage{graphics,latexsym,amssymb,amsmath,amscd}
\usepackage{ifthen,srcltx,amsopn,amssymb,amsfonts}
\usepackage{amssymb,color}
\usepackage{colordvi,fancyhdr}
\usepackage{asymptote}
\usepackage{amscd}
\usepackage[english]{babel}

\usepackage[T1]{fontenc}
\usepackage[utf8x]{inputenc}
\usepackage{psfrag}
\usepackage{delarray,graphicx}
\usepackage{comment}
\usepackage{tikz}
\usepackage[all,cmtip]{xy}
\usepackage{epsfig,graphics}
\usepackage{graphicx}

\usepackage{amssymb,latexsym}
\usepackage{color}
\usepackage{eucal}

\def\[#1\]{\begin{eqnarray*}#1\end{eqnarray*}}

\newcommand{\Z}{\mathbb{Z}}

\DeclareFontFamily{OT1}{rsfs}{}

\DeclareFontShape{OT1}{rsfs}{n}{it}{<-> rsfs10}{}
\DeclareMathAlphabet{\mathscr}{OT1}{rsfs}{n}{it}
\newcommand{\Q}{\mathbb{Q}}

\newcommand{\R}{\mathbb{R}}

\newtheorem{theorem}{Theorem}[section] % 1st argument is your name for it
\newtheorem{lemma}[theorem]{Lemma}     % 2nd argument is what is printed
\newtheorem{corollary}[theorem]{Corollary}
\newtheorem{proposition}[theorem]{Proposition}

\theoremstyle{definition}
\newtheorem{definition}[subsection]{Definition}

\newtheorem{example}{Example}[section]

\theoremstyle{remark}

\newcommand{\Pf}{{\sc Proof}. }
\newcommand{\EPf}{\hbox{}\hfill$\Box$\vspace{.5cm}}
\newcommand\bZ{\mathbb{Z}}
\newcommand\bC{\mathbb{C}}
\newcommand\bR{\mathbb{R}}

\newcommand\SL{\textrm{SL}}
\newcommand\PGL{\textrm{PGL}}
\newcommand\PU{\textrm{PU}}
\newcommand\Fl{\mathcal{F}l}%variété des drapeaux
\newcommand\Conf{\mathcal{C}onf}
\newcommand\Gc{\mathcal{G}c}
\renewcommand\P{\mathbb{P}}
\numberwithin{equation}{subsection}

\newcommand{\PSL}{\mathrm{PSL}}

\begin{document}

\title{Duality and invariants of representations of fundamental groups of 3-manifolds into PGL$(3,\mathbb{C})$}

\author{Elisha Falbel}
\address{
        Institut de Math\'{e}matiques\\
        Universit\'{e} Pierre et Marie Curie\\
        4, Place Jussieu\\
        F-75252, Paris, France\\
        E-mail:elisha.falbel@imj-prg.fr
}

 \author{ Qingxue Wang }\thanks{Q.Wang was supported by NSFC grants \#11171068 and \#11371092.}
 \address{
        School of Mathematical Sciences \\
        LMNS\\
        Fudan  University\\
        Shanghai, 200433\\
        P.R. China\\
        E-mail:qxwang@fudan.edu.cn}

\maketitle

\begin{abstract}
We determine the explicit transformation under duality of generic configurations of four flags  in $\PGL(3,\bC)$ in cross-ratio coordinates. As an application we prove the invariance under duality of an invariant in the Bloch group obtained from decorated triangulations of 3-manifolds.
\end{abstract}

\section{Introduction}

Representations of fundamental groups of manifolds into $\PGL(3,\bC)$ have been studied from various
points of view. As the bulk of past research concentrated in representations of surface groups, some results have been recently addressed for fundamental groups of 3-manifolds and their representations into $\PGL(n,\bC)$ or $\SL(n,\bC)$ (see \cite{MFP,BFG,GTZ,DGG,G,BBI}). One of the main concerns at this point is to collect a sufficient amount of examples which will give directions for further research.  In particular, it is not an easy matter to find examples of non-trivial representations of fundamental groups of hyperbolic manifolds even for $n=3$, that is, $\PGL(3,\bC)$.  But some special representations have been studied before
as they have values in subgroups of  $\PGL(3,\bC)$, notably $\PU(2,1)$ (see \cite{FW} for references), $\SL(3,\bR)$ and $\PGL(2,\bC)$.

Given a simplicial complex $K$ one can define coordinates in a Zariski open set in the character variety following \cite{BFG} (see also \cite{F1,F3,FW,GGZ,DGG,G}).
The coordinates are obtained through a decoration of the simplices by associating a flag in $\mathbb{CP}^{2}$
to each vertex and imposing that compatibility conditions should hold along faces and edges.  The
decorated simplicial complex is denoted by $(K,z)$, where $z$ stands for a set of complex coordinates.
Using these coordinates one can define an invariant $\beta(K,z)$ (see Definition \ref{def beta} in section 4) in the pre-Bloch group $\mathcal{P}(\bC)$. The 5-term relation expresses the fact that $\beta(K,z)$ is invariant under changes of the triangulation.
In the particular  case of manifolds with cusps, one can start with an ideal triangulation and show that
if we choose a sufficiently thin barycentric sub-division we do obtain the full character variety.

There are two natural involutions on $\PGL(3,\bC)$, namely, complex conjugation and the Cartan involution given by $\Theta(g)=(g^{-1})^T=(g^T)^{-1}$, where for $A\in \PGL(3,\bC)$, $A^T$ denotes its transpose. Those involutions determine two actions
on the representation space of a group into $\PGL(3,\bC)$ and also on the corresponding character variety.

On the configuration space of flags in $\mathbb{CP}^{2}$, there are two natural involutions as well. One is induced by the complex conjugation, the other is induced by duality, that is, by the corresponding flags in the dual vector space. These give rise to two actions on the coordinates and the invariant $\beta(K,z)$. See section $3$ for the details. Moreover, the complex conjugation corresponds to the complex conjugate representation and the duality corresponds to the Cartan involution $\Theta$ on the representation space.

Although complex conjugation is more easily understood, the action of duality on $\beta(K,z)$ is more elusive.  Our goal
in this paper is to examine the action of duality on the invariant $\beta(K,z)$.  We prove that
duality changes the invariant by a boundary term  which we determine explicitly (see Proposition \ref{proposition:main}).  In particular, if the simplicial complex has only cusp boundary,
the invariant is fixed by duality (see Theorem \ref{theorem:main}).  It is interesting to note that if the 3-manifold has boundary, then the difference of the pair of dual invariants, $\beta(K,z)-\beta(K,z^*)$, gives a basic invariant of the decorated boundary surface.
 
Motivation for this work follows from data accumulated in \cite{FKR} and in the CURVE project (see \cite{Cu}).
We have listed unipotent decorated triangulations for cusped hyperbolic manifolds obtained by gluing up to
four tetrahedra.  The representations come in quadruplets consisting of one representation and
its orbit under the group $\bZ_2\times \bZ_2$ generated by conjugation and duality.  We have checked numerically that the volumes of dual representations are the same and that prompted us to guess that the underlying invariants in the Bloch group from which the
volume is computed are also the same. We have also checked numerically equality of the volume under duality for random points in the character variety of the figure eight knot (\cite{FGKRT}), giving us
further evidence that duality preserves the invariant $\beta(K,z)$ even in the pre-Bloch group $\mathcal{P}(\bC)$.

Another motivation for this work is the result related to spherical CR structures proved in \cite{FW} (see \cite{J} as a general reference for CR structures).  We considered ideal triangulations of cusped manifolds with a decoration implying  that each simplex can be realized as having the four vertices in the standard sphere in $\bC^2$ (the figure eight knot example is treated in \cite{F2}) and proved  that the volume of the decoration (defined as the the Bloch-Wigner dilogarithm applied to an element $\beta(K,z)$ of the Bloch group which is canonically associated to the decoration) is zero.
Here we obtain a refined result (see Theorem \ref{cr case}).

A third motivation to this work, as pointed out by the referee,  is to understand the Borel class in $H_{c}^3( \PGL(3,\bC), \bR)$ (this is a generator of the continous cohomology group, see \cite{BBI}).  Our result shows that there are two natural dual cocycle representatives, which are different, and we compute the exact coboundary term relating them. We believe that an analogous invariance under duality will also hold for $\PGL(n,\bC)$, $n\geq 4$.

The organization of the paper is the following. In section $2$, we recall the definitions and basic properties of (pre-)Bloch group, dilogarithm function and duality on the configuration space of flags in a projective space. In section $3$, we first review the coordinates defined in \cite{BFG} which parametrize the generic configurations of $3$ and $4$ ordered flags in $\mathbb{CP}^{2}$. Then we study in detail how these coordinates change under duality. The main result in this section is Proposition \ref{proposition:main} which shows that for $T$ a tetrahedron of flags, the difference of the $\mathcal{P}(\bC)$-valued invariant $\beta(T)$ and its dual is a boundary term.  In section $4$, for a decorated simplicial complex $(K,z)$, we study its $\mathcal{P}(\bC)$-valued invariant $\beta(K,z)$ and prove the main Theorem \ref{theorem:main} saying that if the simplicial complex has only cusp boundary, then $\beta(K,z)$ is invariant under duality. Finally in section $5$, we give two applications. One is about spherical CR structures which refines a result in \cite{FW}. The other is that for a cusped hyperbolic manifold $M$ of finite volume, the $\mathcal{P}(\bC)$-valued invariant $\beta(\rho)$ (see Definition \ref{rho}) of a decorated representation $\rho: \pi_1(M)\rightarrow \PGL(3,\bC)$ is invariant under duality.

\section*{Acknowledgements} We thank Antonin Guilloux for discussions on duality and the result in the paper. We are very grateful to an anonymous referee who read our paper very carefully and suggested the beautiful geometric insight on the proof of the Proposition \ref{*prop} part (2). The first author thanks the Department of Mathematics and the Key Laboratory of Mathematics for Nonlinear Sciences (MOE), Fudan University for the financial support under the Senior Visiting Scholarship Project. Both authors thank ANR SGT (Structures G\'eom\'etriques triangul\'ees) for financial support.

\section{Preliminaries}
\subsection{The pre-Bloch group, Bloch group and Dilogarithm}
There are several definitions of Bloch group in the literature and we will mainly follow \cite{S1} here.
\begin{definition}
Let $F$ be a field. The pre-Bloch group $\mathcal{P}(F)$ is the quotient of the free
abelian group $\Z[F\setminus\{0,1\}]$ by the subgroup generated by
the 5-term relation:
\begin{equation}\label{5term}
 [x]-[y]+[\frac{y}{x}]-[\frac{1-x^{-1}}{1-y^{-1}}]+[
  \frac{1-x}{1-y}],\ \  \forall\, x,y\in F\setminus\{0,1\}, x\ne y.
\end{equation}
\end{definition}

For $z\in F\setminus\{0,1\}$, we will still denote by $[z]$ the element it represents in $\mathcal{P}(F)$. If $F$ is algebraically closed, then we have two more relations in $\mathcal{P}(F)$ (\cite{DS} Lemma 5.11): \\ $\forall \, z\in F\setminus\{0,1\}$,
\begin{equation}\label{inv}
 [z^{-1}]=-[z];
\end{equation}
\begin{equation}\label{min}
[1-z]=-[z].
\end{equation}

Now for $a,b \in \bC \setminus\{0,1\}$ and $ab\ne 1$, if we take $x=a^{-1}$, $y=b$ in the 5-term relation (\ref{5term}), by (\ref{inv}), we get the following identity in $\mathcal{P}(\bC)$:
\begin{equation}\label{prod}
[ab]=[a]+[b]+[\frac{1-a}{1-b^{-1}}]+[\frac{1-b}{1-a^{-1}}].
\end{equation}
We will need these identities in later sections.

Consider the tensor product $F^*\otimes_{_\Z} F^*$, where $F^*=F\setminus \{0\}$ is the
multiplicative group of $F$. Let $T=\langle x\otimes y+y\otimes x\ |\ x,y\in F^{*} \rangle$ be the subgroup of $F^{*}\otimes_{_\Z}F^{*}$ generated
by $x\otimes y+y\otimes x$, where $x,y\in F^{*}$.
\begin{definition}
$\bigwedge^2F^{*}= (F^{*}\otimes_{_\Z}F^{*})/T$. For $x,y\in F^{*}$, we will denote by $x\wedge y$ the image of
$x\otimes y$ in $\bigwedge^2F^{*}$.
\end{definition}
By definition, for $x,y \in F^{*}$, we have $x\wedge y=-y\wedge x$ and $2\, x\wedge x=0$.
But $x\wedge x=0$ is
not necessarily true in $\bigwedge^2F^{*}$.

Define a homomorphism $\delta_0: \Z[F\setminus\{0,1\}]\rightarrow \bigwedge^2F^{*}$ as follows: for each generator $z\in F\setminus\{0,1\}$,
$\delta_0 (z)=z\wedge(1-z)$. Then we can check that $\delta_0(\mbox{5-term relation})=0$ (\cite{S1} Lemma 1.1). Hence $\delta_0$ induces a homomorphism
$$\delta: \mathcal{P}(F) \rightarrow {{\bigwedge}^2}F^{*}.$$

\begin{definition}
The Bloch group $\mathcal{B}(F)$ is the kernel of the homomorphism
$\delta$. It is a subgroup of $\mathcal{P}(F)$.
\end{definition}

When $F=\bC$, $\mathcal{P}(\bC)$ (\cite{DS} Theorem 4.16) and  $\mathcal{B}(\bC)$ (\cite{S1})
are uniquely divisible groups and, in fact, are $\Q$-vector
spaces with infinite dimension (\cite{S1}).
In particular they have no torsion. On the other hand, if $F=\R$, then
there exists torsion.  In particular, for all $x\in \R-\{0,1\}$,
the element $[x]+[1-x]\in \mathcal{B}(\R)$ does not depend on $x$ and
has order six (\cite{S1} Proposition 1.1).

Consider the complex conjugation in $\bC$ and its extension
to an involution:
$$
\tau :\Z[\bC\setminus \{0,1\}]
\rightarrow \Z[\bC\setminus \{0,1\}].
$$  As $\tau$
preserves the 5-term relation (\ref{5term}), it induces an involution on the
pre-Bloch group $\mathcal{P}(\bC)$ which we will also denote by $\tau$. Set $$\mathcal{P}(\bC)^{+}=\{x\in \mathcal{P}(\bC)| \tau(x)=x\}$$ and $$\mathcal{P}(\bC)^{-}=\{x\in \mathcal{P}(\bC)| \tau(x)=-x\}.$$ They are the corresponding $\pm 1$-eigenspaces of $\tau$. Then we see
$$
 \mathcal{P}(\bC)=\mathcal{P}(\bC)^{+}\oplus \mathcal{P}(\bC)^{-}.
$$
Similarly,
$$
 \mathcal{B}(\bC)=\mathcal{B}(\bC)^{+}+\mathcal{B}(\bC)^{-},
$$
where $\mathcal{B}(\bC)^{+}=\{x\in \mathcal{B}(\bC)| \tau(x)=x\}$ and $\mathcal{B}(\bC)^{-}=\{x\in \mathcal{B}(\bC)| \tau(x)=-x\}$.
Clearly $\mathcal{B}(\bC)^{+}$ (resp. $\mathcal{B}(\bC)^{-}$ ) is a subgroup of $\mathcal{P}(\bC)^{+}$ (resp. $\mathcal{P}(\bC)^{-}$).

\begin{definition}\label{BWf} The Bloch-Wigner dilogarithm function $D:\bC-\{0,1\}\rightarrow \bR$ is defined as follows: for $x\in \bC-\{0,1\}$, define
$$
  D(x)=\arg{(1-x)}\log{|x|}-\Im (\int_{0}^{x}\log{(1-t)}\frac{dt}{t}),
$$
where for a complex number $z$, $\Im z$ is its imaginary part.
\end{definition}

It is well-defined and real analytic on $\bC-\{0,1\}$ and
extends to a continuous function on the projective line $\mathbb{CP}^{1}$ by
defining $D(0) = D(1) = D(\infty) = 0$. It is well-known
that it satisfies the 5-term relation (\ref{5term}). Hence it linearly extends to a well-defined homomorphism:
$$
  D: \mathcal{P}(\bC)\rightarrow \R,
$$
given by
$$
 D(\sum_{i=1}^{k}n_i[x_i])=\sum_{i=1}^{k}n_iD(x_i).
$$
Since $D(\bar{z})$=-$D(z)$ for any $z\in \bC$, we obtain that
\begin{equation}\label{dplus}
D(u)=0, \, \forall \, u\in \mathcal{P}(\bC)^{+}.
\end{equation}

For the proofs of the above-mentioned properties of the
dilogarithm function $D$, see \cite{B1} Lecture 6.

\subsection{ Configuration space of flags and duality}
Let $V$ be a complex vector space of dimension $n$. Let $\P(V)$  be the projective space of $V$. We have the natural map $\pi: V\setminus \{0\} \rightarrow \P(V)$ which sends each non-zero vector $v$ to the line $[v]$ generated by $v$. For a subspace $W$ of $V$, we will denote $\pi(W\setminus \{0\})$ by $\P(W)$. \emph{A flag of $\P(V)$} is defined to be an ascending sequence of linear subspaces of $\P(V)$ : $\P(V_1)\subset \P(V_2)\subset \cdots \subset \P(V_n)=\P(V)$, where $\{0\}\subset V_1\subset V_2\subset \cdots \subset V_n=V$ is a (complete) flag of $V$, i.e., $\dim V_i=i$, $1\leq i\leq n$. Let $\Fl(V)$ be the set of all flags of $\P(V)$. For a positive integer $m$, we define $\Fl_m(V)=\{(f_1,\cdots,f_m)|f_i\in \Fl(V), 1\leq i\leq m\}$ to be the set of $m$ ordered flags of $\P(V)$. Since $\PGL(n,\bC)$ acts naturally on $\Fl(V)$, it acts term-by-term on $\Fl_m(V)$. The configuration space of $m$ ordered flags of $\P(V)$ is defined as the orbit space of $\Fl_m(V)$ by the action of $\PGL(n,\bC)$. We will denote it by $\Conf_m$. For $(F_1,\cdots,F_m)\in \Fl_m(V)$, we will denote by $[(F_1,\cdots,F_m)]$ its orbit in $\Conf_m$.

%\subsection{Duality}
Let $V^*$ be the dual vector space of $V$, this is, the space of all linear functionals of $V$. Let $\P(V^*)$ be the projective space of $V^*$. Let $\Fl(V^*)$ be the set of all flags of $\P(V^*)$. Given a subspace $W$ of $V$ with $\dim{W}=l$, we define $W'=\{f \in V^*|f(x)=0,\, \forall x\in W\}$. It is a subspace of $V^*$ with $\dim{W'}=n-l$. Let $\{0\}\subset V_1\subset V_2\subset \cdots \subset V_n=V$ be a flag of $V$. Then we have $\{0\}=V_n'\subset V_{n-1}' \subset \cdots \subset V_1' \subset \{0\}'=V^*$ is a flag of $V^*$.  Now we define a map $\sigma: \Fl(V)\rightarrow \Fl(V^*)$ as follows: for a flag $F =\P(V_1)\subset \P(V_2)\subset \cdots \subset \P(V_n)=\P(V)$,
 $$\sigma(F)=\P(V_{n-1}') \subset \P(V_{n-2}')\cdots \subset \P(V_1')\subset \P(V^*).
 $$
If we fix a basis $\{e_1,\cdots,e_n \}$ of $V$, we have the dual basis $\{e_1^*,\cdots, e_n^*\}$ of $V^*$ and an isomorphism between $V$ and $V^*$ which maps $e_i$ to $e_i^*$, $1\leq i\leq n$. Under this isomorphism, we identify $\Fl(V^*)$ with $\Fl(V)$. Now the above map $\sigma: \Fl(V)\rightarrow \Fl(V^*)=\Fl(V)$ extends to an involution $\Theta_m: \Conf_m \rightarrow \Conf_m$. That is, for $[(F_1,\cdots,F_m)]\in \Conf_m$, $$\Theta_m([(F_1,\cdots,F_m)])=[(\sigma(F_1),\cdots,\sigma(F_m))].$$
It is straightforward to check that $\Theta_m$ is independent of the chosen basis and $\Theta_m\circ \Theta_m=id$.

\begin{definition}\label{d conf}
We call $\Theta_m$ \emph{the duality} on the space of configurations of $m$ ordered flags of $\P(V)$. For $C=[(F_1,\cdots,F_m)]\in \Conf_m$, we call $\Theta_m(C)=[(\sigma(F_1),\cdots,\sigma(F_m))]$ its dual configuration, denoted by $C^*$.
\end{definition}
\section{Coordinates for Configuration space of flags and duality in $\mathbb{CP}^{2}$}\label{ss:coord-para}
From now on, we will only consider the case $V=\bC^3$. We will fix the standard basis $\{e_1=(1,0,0), e_2=(0,1,0), e_3=(0,0,1)\}$ of $V$ and the corresponding dual basis $\{e_1^*, e_2^*,e_3^*\}$. Then $\P(V)=\mathbb{CP}^{2}$ is the projective plane.

\subsection{The action of $\PGL(3,\bC)$ and duality}
Note that a hyperplane of $V$ corresponds to the kernel of a linear functional, so we obtain
$$\Fl(V)=\{([x],[f])\in \P(V)\times \P(V^*) \, | \, f(x)=0\}.$$ Here $([x],[f])$ with $f(x)=0$ represents the flag $\{[x]\}\subset \P(\ker{f})$ of $\P(V)$. Hence the natural action of $\PGL(3,\bC)$ on $\Fl(V)$ has the following form. For $A\in \PGL(3,\bC)$ and $([x],[f])\in \Fl(V)$ with projective coordinates
$x=[a_0,a_1,a_2]$, $f=[u_0,u_1,u_2]$, $$A([x],[f])=([\tilde{a}_0, \tilde{a}_1,\widetilde{a}_2],[\tilde{u}_0,\tilde{u}_1,\tilde{u}_2]),$$ where $$
(\tilde{a}_0, \tilde{a}_1,\widetilde{a}_2)^T=A\cdot (a_0,a_1,a_2)^T, \,\,  (\tilde{u}_0,\tilde{u}_1,\tilde{u}_2)^T=(A^T)^{-1}\cdot (u_0,u_1,u_2)^T.$$
Here $^T$ means the transpose of a matrix.

Under the canonical isomorphism between $V$ and its double dual $V^{**}$, we get
$$\Fl(V^*)=\{([g],[y])\in \P(V^*)\times \P(V) \, | \, g(y)=0\}.$$ By definition, the map $\sigma: \Fl(V)\rightarrow \Fl(V^*)=\Fl(V)$ is simply switching the point and the line, i.e. $\forall F=([x],[f])\in \Fl(V)$,
$$\sigma(F)=([f],[x]).$$
Hence, for $C=[([x_1],[f_1]),\cdots, ([x_m],[f_m])]\in \Conf_m$, the duality $\Theta_m:\Conf_m\rightarrow \Conf_m$ is given by
$$ C^*=\Theta_m(C)=[([f_1],[x_1]),\cdots, ([f_m],[x_m])].$$

For a complex number $z$, we denote by $\bar{z}$ its complex conjugation. Now we have another involution on $\Conf_m$ induced by complex conjugation.  Namely for $C=[([x_1],[f_1]),\cdots, ([x_m],[f_m])]\in \Conf_m$, we define its complex conjugation $$\bar{C}=[([\bar{x}_1],[\bar{f}_1]),\cdots, ([\bar{x}_m],[\bar{f}_m])],$$ where for $([x],[f])\in \Fl(V)$ with projective coordinates $x=[a_0,a_1,a_2]$, $f=[u_0,u_1,u_2]$, we set $\bar{x}=[\bar{a}_0,\bar{a}_1,\bar{a}_2]$, $f=[\bar{u}_0,\bar{u}_1,\bar{u}_2]$. Clearly the map $\iota: \Conf_m \rightarrow \Conf_m$, $\iota(C)=\bar{C}$ is an involution for any $m\geq 1$. By definition, $\iota$ and $\Theta_m$ commute, i.e., for $C\in \Conf_m$, we have
$$ \Theta_m(\iota(C))=\iota(\Theta_m(C)).$$

Clearly $\Conf_1$ consist of one point. Next, we will review the coordinates defined in \cite{BFG} which parametrize the generic configurations in $\Conf_2$, $\Conf_3$ and $\Conf_4$. Then we will study the changes of these coordinates under the duality $\Theta_3$ and $\Theta_4$.

\begin{definition}
Let $[([x_1],[f_1]),\cdots, ([x_m],[f_m])]\in \Conf_m$. We say it is \emph{a generic configuration} if $f_i(x_j)\ne 0$, $\forall i\ne j$ and $[x_1],\cdots,[x_m]$ are in generic position in $\P(V)$ (i.e. no three of them lie on a line in $\mathbb{CP}^{2}$). We say it is \emph{a very generic configuration} if it is generic and the lines $[f_1],\cdots,[f_m]$ are in generic position as well. We will denote by $\Gc_m$ the configuration space of $m$ ordered generic flags of $\P(V)$.
\end{definition}

Note that the configuration space of $m$ very generic flags is a Zariski open subset of $\Gc_m$. The image of a generic but not very generic configuration under $\Theta_m$ might be no longer generic. Thus, when restricting to $\Gc_m$, $\Theta_m$ is a birational automorphism of order dividing $2$.

\subsection{Coordinates for $\Gc_3$ and duality}
Given a generic configuration $C=[([x_1],[f_1]),([x_2],[f_2]), ([x_3],[f_3])]\in \Gc_3$, we have its triple ratio
$$
z_{123}=\frac{f_1(x_2)f_2(x_3)f_3(x_1)}{f_1(x_3)f_2(x_1)f_3(x_2)}\in \bC^{*}=\bC\setminus \{0\},
$$
which classifies the configuration space of three ordered generic flags of $\P(V)$. Indeed, since the points $[x_1]$,$[x_2]$,$[x_3]$ and $\P(\ker{f_1}\cap \ker{f_2})$ are in general position, there is a unique projective transformation mapping them to $[1,0,0]$,$[0,1,0]$,$[0,0,1]$ and $[1,1,1]$ respectively. Therefore the flag becomes
$$C= [([1,0,0],[0,1,-1]),([0,1,0],[1,0,-1]),([0,0,1],[z,1,0])],$$
where $z=z_{123}$ is the triple ratio. By definition, $C$ is not very generic if and only if $[f_1]$,$[f_2]$,$[f_3]$ are collinear, that is, they pass through a common point. In coordinates, this is the equation $0=\det{(f_1,f_2,f_3)}=-1-z$. Thus there is only one generic but not very generic configuration $C_0 \in \Conf_3$, whose triple ratio $z=-1$.

Suppose $C=[([x_1],[f_1]),([x_2],[f_2]),([x_3],[f_3])]$ satisfies the condition $f_i(x_j)\neq 0$ for $i\neq j$.   The triple ratio is well defined in this case and $\Theta_3(C)=[([f_1],[x_1]),([f_2],[x_2]), ([f_3],[x_3])]$ satisfies the same property.
Let $\sigma$ be a permutation of $\{1,2,3\}$. We see that if $\sigma$ is  an even permutation, then $\displaystyle{z_{\sigma(1)\sigma(2)\sigma(3)}=z_{123}}$; if it is odd, then $\displaystyle{z_{\sigma(1)\sigma(2)\sigma(3)}=\frac{1}{z_{123}}}$.

The triple ratio of $\Theta_3(C)$ is equal to
$$
z^*_{123}=\frac{f_2(x_1)f_3(x_2)f_1(x_3)}{f_3(x_1)f_1(x_2)f_2(x_3)}=\frac{1}{z_{123}}.
$$

For $C\in \Gc_3$, we see from the definition the triple ratio of its complex conjugation $\bar{C}$ is equal to the complex conjugation of the triple ratio of $C$. Summarizing the above discussions, we have the following
\begin{proposition}\label{3flags}
(1). $\Gc_3$ is identified with $\bC^*$ by the triple ratio. $C\in \Gc_3$ is very generic if and only if its triple ratio in addition satisfies $z\ne -1$. \\
(2). $\Theta_3$ extends to an involution on $\Gc_3$. Moreover by the triple ratio, it corresponds to the involution $t:\bC^*\rightarrow \bC^*$ given by $$t(z)=\frac{1}{z}, \,\, \forall z\in \bC^* .$$
(3). By triple ratio, $\iota:\Gc_3 \rightarrow \Gc_3$ corresponds to the complex conjugation $$c:\bC^*\rightarrow \bC^*,\;\;c(z)=\bar{z}.$$
\end{proposition}

\subsection{Coordinates for $\Gc_4$ and duality}
Following \cite{BFG}, we call a (very) generic configuration in $\Gc_4$ a (very) generic tetrahedron of flags. Figure 1 displays the coordinates on a tetrahedron of flags.

Let $T=([x_i],[f_i])_{1\leq i \leq 4}$ be a generic tetrahedron of flags. There are $12$ coordinates for the $12$ oriented edges of the tetrahedron and $4$ coordinates for the $4$ oriented faces of the tetrahedron. Let's recall the definition of the edge coordinate $z_{ij}(T)$ associated to the edge $ij$, $1\leq i\ne j\leq 4$. We first choose $k$ and $l$ such that the permutation $(1,2,3,4) \mapsto (i,j,k,l)$ is even. Then the set of all the lines in $\P_2(\bC)$ passing through the point $x_i$ is a projective line. We have four distinct points on this projective line: the line $\ker{f_i}$ and the three lines $x_ix_l$ passing through $x_i$ and $x_l$ for $l\neq i$. We define $z_{ij}$ as the cross-ratio of these four points
by
$$z_{ij}(T) := X(\ker{f_i},x_ix_j,x_ix_k,x_ix_l),$$
where we define the cross-ratio $X(x_1,x_2,x_3,x_4)$ of four distinct points
$x_1$,$x_2$,$x_3$,$x_4$ on a projective line to be the value at $x_4$ of the linear fractional transformation mapping
$x_1$ to $\infty$, $x_2$ to $0$, and $x_3$ to $1$. Hence we have the formula
$$X(x_1,x_2,x_3,x_4)=\frac{(x_1-x_3)(x_2-x_4)}{(x_1-x_4)(x_2-x_3)}. $$
By \cite{BFG} Lemma 2.3.1, we have another description for $z_{ij}(T)$:
\begin{equation}\label{z coord}
z_{ij}(T)=\frac{f_i(x_k)\det{(x_i,x_j,x_l)}}{f_i(x_l)\det{(x_i,x_j,x_k)}}.
\end{equation}

For each face $(ijk)$ oriented as the boundary of the tetrahedron $(1234)$, we define the face coordinate $z_{ijk}(T)$ to be the triple ratio of $([x_i],[f_i]),([x_j],[f_j]),([x_k],[f_k])$. That is,
$$
z_{ijk}(T) = \frac{f_i(x_j)f_j(x_k)f_k(x_i)}{f_i(x_k)f_j(x_i)f_k(x_j)}.
$$

By (\ref{z coord}), we obtain

\begin{equation}\label{faceedge}
z_{ijk}(T) =- z_{il}(T)z_{jl}(T)z_{kl}(T).
\end{equation}
This means $z_{ijk}(T)$ is the opposite of the product of all edge coordinates ``leaving'' this face.

When gluing tetrahedra along a face it is important to notice that if the same face $(ikj)$ (with opposite orientation) is common to a second tetrahedron $T'$ then $$z_{ikj}(T') = \frac{1}{z_{ijk}(T)}.$$

\begin{figure}\label{tetra}
  \centering
{\scalebox{1}{\includegraphics{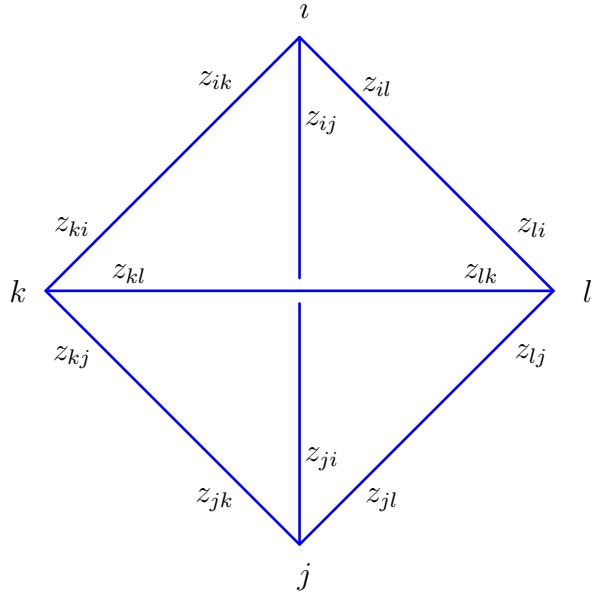}}}
  \caption{The $z$-coordinates.}
    \label{eight}
\end{figure}

By the basic properties of cross-ratios, we see that the edge coordinates leaving a vertex $i$ are related by:
\begin{equation} \label{Cros}
\begin{split}
& z_{ik}(T) =  \frac{1}{1-z_{ij}(T)},\\
& z_{il}(T) =   1-\frac{1}{z_{ij}(T)}.
\end{split}
\end{equation}
Observe that the relations follow a cyclic order around each vertex which is defined by the orientation of the tetrahedron. Moreover the above relation implies that at each vertex $i$, the edge coordinates leaving $i$ satisfy:
\begin{equation} \label{e prod}
 z_{ij}(T)z_{ik}(T)z_{il}(T)=-1.
\end{equation}

\begin{proposition}\label{4flags}
(1). A generic tetrahedron of flags $T\in \Gc_4$ is uniquely determined by the $4$ edge coordinates $(z_{12}(T),z_{21}(T),z_{34}(T),z_{43}(T))$ which are in $\bC \setminus \{0,1\}$. Hence $$\Gc_4 \simeq  (\bC \setminus \{0,1\})^4.$$\\
(2). $T\in \Gc_4$ is very generic if and only if the face invariants satisfy $$z_{123}(T)\ne -1,\;z_{124}(T)\ne -1,\;z_{134}(T)\ne -1,\;z_{234}(T)\ne -1.$$
\end{proposition}

\Pf (1). Note that since $T$ is generic, each edge coordinate $z_{ij}(T)\in \bC \setminus \{0,1\}$. By \cite{BFG} Proposition 2.4.1, under the action of $\PGL(3,\bC)$, $T$ is uniquely given by the following four flags $([x_i],[f_i])$, ${1\leq i \leq 4}$ with: $[x_1]=[1,0,0],\; [f_1]=[0,1-\frac{1}{z_{12}(T)},-1]$; $[x_2]=[0,1,0],\; [f_2]=[1-z_{21}(T),0,-1]$; $[x_3]=[0,0,1],\; [f_3]=[z_{34}(T),-1,0]$; $[x_4]=[1,1,1],\; [f_4]=[\frac{1}{1-z_{43}(T)},1-\frac{1}{1-z_{43}(T)},-1]$.\\

(2). By definition, $T=([x_i],[f_i])_{1\leq i \leq 4}$ is very generic if and if only any three of the $[f_i]$, $1\leq i\leq 4$ are not collinear, i.e. any three flags in $T$ form a very generic configuration in $\Gc_3$. Then it follows from Proposition \ref{3flags} and the definition of the face coordinates.
\EPf

Note that the above isomorphism $\Gc_4 \simeq  (\bC \setminus \{0,1\})^4$ is not canonical. In fact one can choose one edge coordinate at each vertex, for instance, $(z_{13}(T),z_{21}(T),z_{32}(T),z_{43}(T))$ will also give an isomorphism.

Let $T\in \Gc_4$ be a generic tetrahedron with edge coordinates $z_{ij}(T)$, $1\leq i\ne j\leq 4$. We define a $\mathcal{P}(\bC)$-valued invariant:

\begin{definition}\label{t inv}
$\beta(T)=[z_{12}(T)]+[z_{21}(T)]+[z_{34}(T)]+[z_{43}(T)]\in \mathcal{P}(\bC)$. We define the volume of $T$ to be $\frac{1}{4}D(\beta(T))$, where $D$ is the Bloch-Wigner dilogarithm function in Definition \ref{BWf}. We denote it by $Vol(T)$.
\end{definition}

Let $F=([x],[f])\in \Fl(V)$ be a flag. We set $F^*=\sigma(F)=([f],[x])$. Let $T=(F_i)_{1\leq i\leq 4}\in \Gc_4$ be a very generic tetrahedron of flags. Then its dual $T^*=\Theta_4(T)=(F_i^*)_{1\leq i\leq 4}$ is also very generic. Next we study the relations between the coordinates of $T$ and $T^*$.

\begin{proposition}\label{*prop}
Let $T=(F_i)_{1\leq i\leq 4}$ be a very generic tetrahedron of flags with edge coordinates $z_{ij}$ and face coordinates $z_{ijk}$. Let $T^*$ be its dual configuration with edge coordinates $z^*_{ij}$ and face coordinates $z^*_{ijk}$.  Then\\
(1). For the face coordinates, we have $\displaystyle{z^*_{ijk}=\frac{1}{z_{ijk}}}.$\\
(2). For the edge coordinates, choose $k$,$l$ such that the permutation $(1,2,3,4) \mapsto (i,j,k,l)$ is even, then we have
$$
 z^*_{ij} =z_{ji}\frac{1-z_{ik}z_{jk}z_{lk}}{1-\frac{1}{z_{il}z_{jl}z_{kl}}}=z_{ji}\frac{1+z_{jil}}{1+\frac{1}{z_{ijk}}}.
$$
\end{proposition}

\Pf (1). This follows from Proposition \ref{3flags}. \\
(2). Since the second equality follows from (\ref{faceedge}), we need to show the first equality. We will compute explicitly the dual coordinates and obtain the formula after some manipulations. Indeed, recall that in the proof of Proposition \ref{4flags}, we have the following coordinates for $T$:
$$
F_1= ([1,0,0], [0,z_1,-1]), F_2= ([0,1,0], [z_2,0,-1]),
$$
$$
F_3= ([0,0,1], [z_3,0,-1]), F_4= ([1,1,1], [z_4,1-z_4,-1]),
$$
where
$$
z_1=1-\frac{1}{z_{12}};\, z_2=1-z_{21};\, z_3=z_{34};\, z_4=\frac{1}{1-z_{43}}.
$$
To get the corresponding coordinates for $z^*_{ij}$, we do the following. First there is a unique $A\in \PGL(3,\bC)$ which maps $[0,z_1,-1]$ to $[1,0,0]$, $[z_2,0,-1]$ to $[0,1,0]$, $[z_3,0,-1]$ to $[0,0,1]$ and  $[z_4,1-z_4,-1]$ to $[1,1,1]$. A direct calculation shows that
$$
(A^{-1})^{T}=
\left(
  \begin{array}{ccc}
    0 & az_1 & -a \\
    bz_2 & 0 & -b \\
    cz_3 & -c & 0 \\
  \end{array}
\right),
$$
where
$$
a=\frac{z_4+z_3(1-z_4)-z_2}{z_1z_3-z-2},\; b=\frac{-z_4-z-3(1-z_4)+z_1}{z_1z_3-z-2},\; c=\frac{z_1z_4+z_2(1-z_4)-z_1z_2}{z_1z_3-z-2}.
$$
Now we will show how to find the formula of $z^*_{12}$. The other cases are similar and we will omit the details. Apply $(A^{-1})^{T}$ to $(1,0,0)^T$, we obtain
$$
z^*_{12}=\frac{cz_3}{cz_3+bz_2}=z_{34}\frac{z_{12}z_{21}+z_{21}z_{43}-z_{21}-z_{43}}
 {z_{12}z_{21}+z_{12}z_{34}-z_{12}-z_{34}}=z_{34}\frac{z_{21}(z_{12}-1)+z_{43}(z_{21}-1)}
 {z_{12}(z_{21}-1)+z_{34}(z_{12}-1)}.
$$
By (\ref{Cros}), $z_{12}-1=-1/z_{13}$, $z_{21}-1=-1/z_{24}$, hence
$$
z^*_{12}=z_{34}\frac{z_{21}z_{24}+z_{43}z_{13}}
 {z_{12}z_{13}+z_{34}z_{24}}.
$$
Now by (\ref{e prod}), we have $z_{12}z_{13}z_{14}=-1$, $z_{21}z_{23}z_{24}=-1$. Therefore,
\begin{alignat}{2}
z^*_{12}&=z_{34}\frac{z_{21}z_{24}+z_{43}z_{13}}{z_{12}z_{13}+z_{34}z_{24}}=z_{34}\frac{z_{14}(z_{21}z_{24}+z_{43}z_{13})}
 {z_{14}(z_{12}z_{13}+z_{34}z_{24})}\notag\\
 &=z_{34}\frac{z_{14}(z_{21}z_{24}+z_{43}z_{13})}{-1+z_{14}z_{24}z_{34}}=\frac{z_{21}+\frac{z_{43}z_{13}z_{21}z_{23}}{z_{24}z_{21}z_{23}}}{1-\frac{1}{z_{14}z_{24}z_{34}}}\notag \\
 &=z_{21}\frac{1-z_{13}z_{23}z_{43}}{1-\frac{1}{z_{14}z_{24}z_{34}}}\notag.
\end{alignat}
\EPf

% REMARK DELETED

Now from Proposition \ref{*prop}, it follows immediately
\begin{corollary}\label{cor:ij=kl}
 Let $T=(F_i)_{1\leq i\leq 4}$ be a very generic tetrahedron of flags with edge coordinates $z_{ij}$. Let $T^*$ be its dual configuration with edge coordinates $z^*_{ij}$.  Then for $\{i,j,k,l\}=\{1,2,3,4\}$
 $$
 z^*_{ij}z^*_{ji} = z_{kl} z_{lk}.\\
$$
 \end{corollary}

\begin{figure}\label{configuration}
  \centering
{\scalebox{1}{\includegraphics{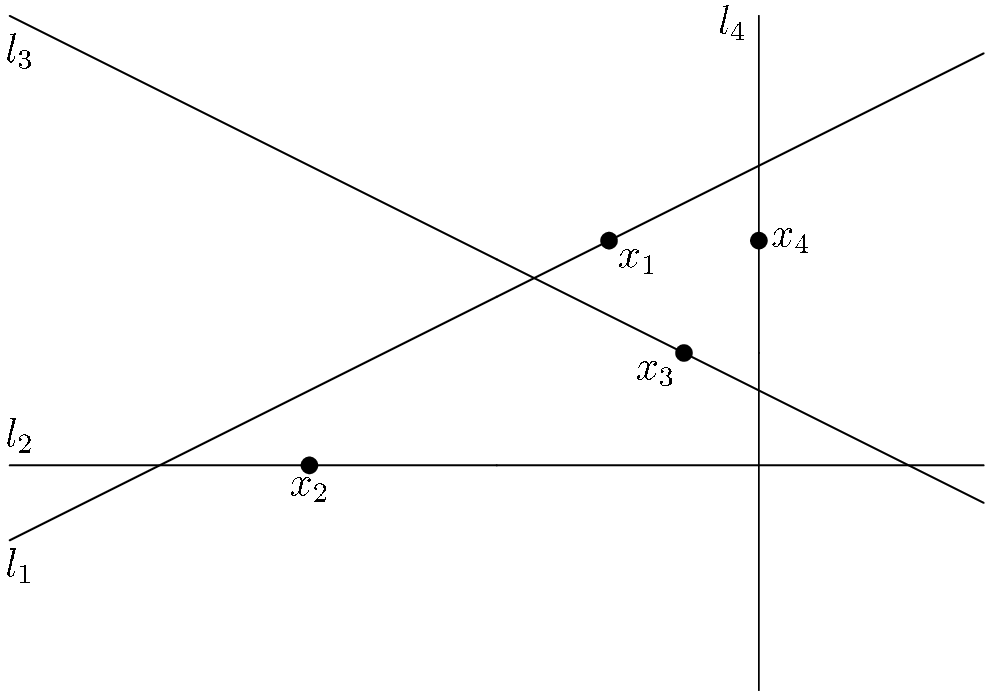}}}
  \caption{The configuration of four flags.  Here $l_i=\ker (f_i)$.}
  \end{figure}

  \begin{figure}\label{z12}
  \centering
{\scalebox{1}{\includegraphics{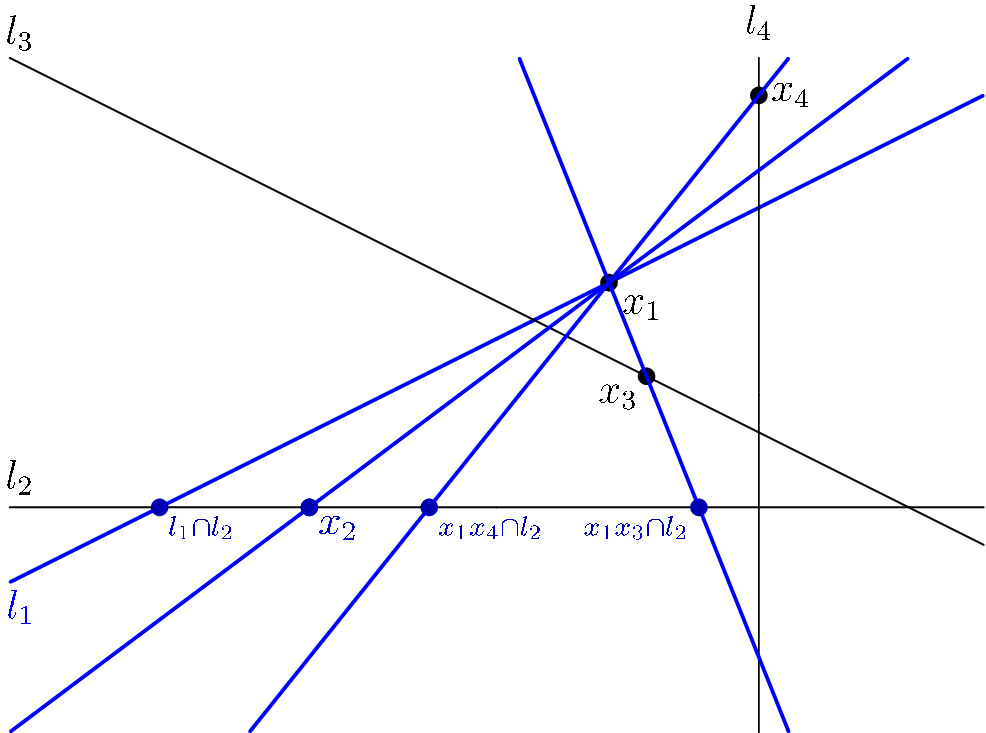}}}
 \caption{The four lines involved in the definition of the cross-ratio $z_{12}=X(l_1,x_1x_2,x_1x_3,x_1x_4)$.  It can be computed as  the cross ratio of the intersection of these lines with the line $l_2$.}
 \end{figure}

  \begin{figure}\label{z21*}
  \centering
{\scalebox{1}{\includegraphics{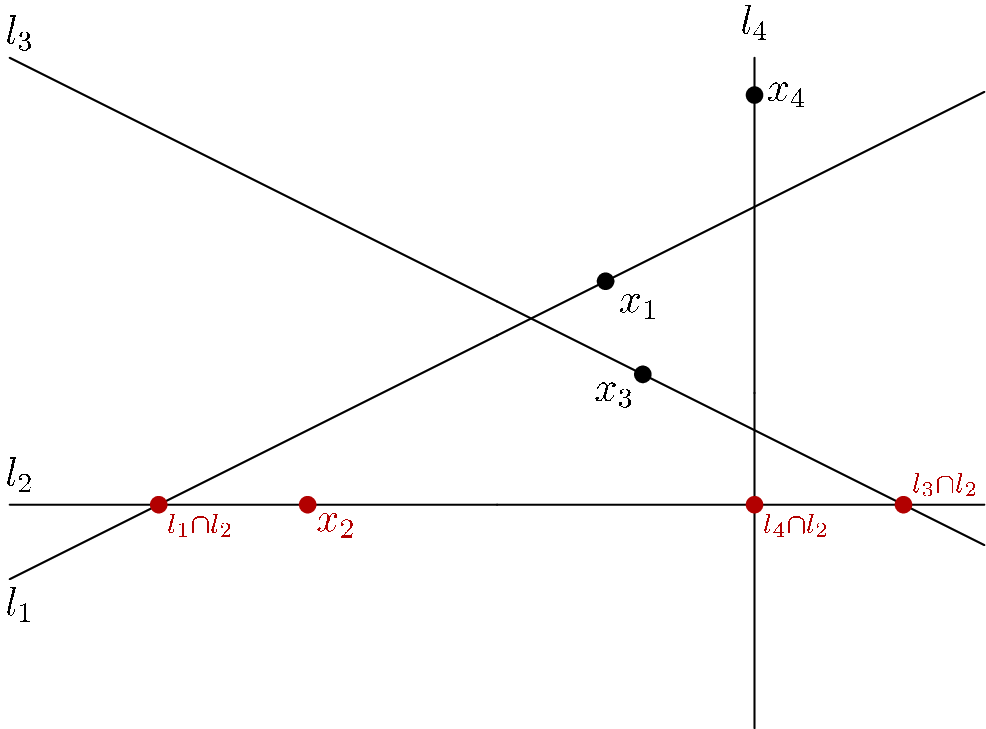}}}
\caption{The four points in the line $l_2$ involved in the definition of the cross-ratio $z^*_{21}=X(x_2,l_2\cap l_1,l_2\cap l_4,l_2\cap l_3)$.}
 \end{figure}

The following proof of Proposition \ref{*prop} part (2) was suggested to us by an anonymous referee.  It makes computations easier and more transparent and geometric. The main observation is that the triple ratio is also a cross-ratio (cf. Lemma 3.8 \cite{Go}) :

\begin{lemma}\label{triple-cross} Let $([x_{i}],[f_i])_{1\leq i\leq 3}$ be a generic triple of flags. Set $l_i=\ker{f_i}$, $1\leq i\leq 3$. On the line $l_2$, there are 4 distinct points $l_2\cap l_3$, $l_2\cap l_1$, $l_2\cap x_1x_3$, $x_2$ .  Then
$$
-z_{123}=X(l_1\cap l_2,l_2\cap l_3,x_2, l_2\cap x_1x_3), \; \; 1+z_{123}=X(l_1\cap l_2, x_2, l_3\cap l_2, l_2\cap x_1x_3).
$$
\end{lemma}

\Pf
Since both the triple ratio and the cross-ratio are invariant under the action of $\PGL(3,\bC)$, from the discussion in Section 3.2,  we can take the 3 flags as follows:
$$ C= [([1,0,0],[0,1,-1]),([0,1,0],[1,0,-1]),([0,0,1],[z,1,0])].$$
Then $-z_{123}=-z$ and the 1st identity follows from the direct computation of the right-hand side. For 4 distinct points $a$,$b$,$c$,$d$ on a projective line, we have $X(a,b,c,d)=1-X(a,c,b,d)$. Now the 2nd identity follows from the 1st one.
\EPf

\Pf(of Proposition \ref{*prop} part (2))
Set $l_i=\ker{f_i}$, $1\leq i\leq 4$. Let us prove the formula
$$
 z^*_{21} =z_{12}\frac{1+z_{123}}{1+\frac{1}{z_{214}}}=z_{12}\frac{1+z_{123}}{1+z_{124}}.
$$
 
By Lemma \ref{triple-cross}, we have
$$
1+z_{123}=X(l_1\cap l_2, x_2, l_3\cap l_2, l_2\cap x_1x_3),
$$

$$
\frac{1}{1+z_{124}}=\frac{1}{X(l_1\cap l_2, x_2, l_4\cap l_2, l_2\cap x_1x_4)}=X(l_1\cap l_2, x_2, l_2\cap x_1x_4,l_4\cap l_2 ).
$$
By the definition (See Figure 3 and Figure 4), we have
$$
z_{12}=X(l_1\cap l_2, x_2,l_2\cap x_1x_3, l_2\cap x_1x_4),
$$
$$
 z^*_{21}=X(x_2,l_1\cap l_2,l_4\cap l_2,l_3\cap l_2)=X(l_1\cap l_2,x_2,l_3\cap l_2,l_4\cap l_2).
$$

Notice that now all the cross-ratios are from the six distinct points on the line $l_2$, namely,
$$ a_1=l_1\cap l_2,\; a_2=x_2,\; b_1=l_3\cap l_2,\; b_2=l_2\cap x_1x_3,\;b_3=l_2\cap x_1x_4,\; b_4=l_4\cap l_2.
$$

The formula follows from the following identity for cross-ratios.
$$
X(a_1,a_2,b_1,b_4)=X(a_1,a_2,b_1,b_2)X(a_1,a_2,b_2,b_3)X(a_1,a_2,b_3,b_4)
$$
The other cases are similar.
\EPf

 Now we can prove the following main proposition which compares $\beta(T)$ and $\beta(T^*)$ in $\mathcal{P}(\bC)$.

\begin{proposition} \label{proposition:main} Let $T=(F_i)_{1\leq i\leq 4}$ be a very generic tetrahedron of flags with edge coordinates $z_{ij}$ and face coordinates $z_{ijk}$. Let $T^*$ be its dual configuration. Then in ${\mathcal{P}}(\bC)$, we have
$$
\beta(T)=\beta(T^*)+[-z_{123}]+[-z_{243}]+[-z_{134}]+[-z_{142}].
$$
\end{proposition}
 \Pf
There are two cases. \\
Case (I): $z_{12}z_{21}\ne 1$ and $z_{34}z_{43}\ne 1$.
By (\ref{prod}), $\forall\, x,y\in \bC \setminus\{0,1\}$ and $xy\ne 1$, we have

\begin{equation}\label{prop1}
 [xy]=[x]+[y]+[\frac{1-x}{1-y^{-1}}]+[\frac{1-y}{1-x^{-1}}].
\end{equation}

Therefore we obtain
$$
[z_{12}z_{21}]=[z_{12}]+[z_{21}]+[\frac{1-z_{12}}{1-1/z_{21}}]+[\frac{1-z_{21}}{1-1/z_{12}}]
 $$
 which by (\ref{Cros}) simplifies to
 $$
[z_{12}]+[z_{21}]= [z_{12}z_{21}]+[z_{13}z_{23}]+[z_{24}z_{14}].
 $$
Analogously, we have
 $$
[z_{34}]+[z_{43}]= [z_{34}z_{43}]+[z_{31}z_{41}]+[z_{32}z_{42}].
 $$
From the expressions above and Corollary
 \ref{cor:ij=kl} we obtain that
 $$\beta(T)-\beta(T^*)=[z_{13}z_{23}]+[z_{24}z_{14}]+[z_{31}z_{41}]+[z_{32}z_{42}]-
 ([z^*_{13}z^*_{23}]+[z^*_{24}z^*_{14}]+[z^*_{31}z^*_{41}]+[z^*_{32}z^*_{42}]).
$$
 Indeed, Corollary \ref{cor:ij=kl} implies that all terms of the form $[z_{ij}z_{ji}]$ cancel out.

 Now we observe that, by (\ref{prop1}) and (\ref{Cros}) again,
 \begin{alignat}{2}
 [z_{13}z_{23}]=[z_{13}z_{23}z_{43}\frac{1}{z_{43}}] &=[z_{13}z_{23}z_{43}]+[\frac{1}{z_{43}}]+
 [\frac{1-z_{13}z_{23}z_{43}}{1-z_{43}}]+[\frac{1-1/z_{43}}{1-1/z_{13}z_{23}z_{43}}]\notag \\
 &=[z_{13}z_{23}z_{43}]+[\frac{1}{z_{43}}]+[z_{42}(1-z_{13}z_{23}z_{43})]+[z_{41}\frac{z_{13}z_{23}z_{43}}{z_{13}z_{23}z_{43}-1}].\notag
 \end{alignat}
By (\ref{faceedge}), we have $z_{13}z_{23}z_{43}=-z_{142}$. Hence we can rewrite $[z_{13}z_{23}]$ as
$$
 [z_{13}z_{23}]=[-z_{142}]+[\frac{1}{z_{43}}]+[z_{42}(1+z_{142})]+[z_{41}\frac{z_{142}}{(1+z_{142})}]
$$
Analogously we have
$$[z_{14}z_{24}]=[-z_{123}]+[\frac{1}{z_{34}}]+[z_{31}(1+z_{123})]+[z_{32}\frac{z_{123}}{(1+z_{123})}]$$
$$[z_{31}z_{41}]=[-z_{243}]+[\frac{1}{z_{21}}]+[z_{24}(1+z_{243})]+[z_{23}\frac{z_{243}}{(1+z_{243})}]$$
$$[z_{32}z_{42}]=[-z_{134}]+[\frac{1}{z_{12}}]+[z_{13}(1+z_{134})]+[z_{14}\frac{z_{134}}{(1+z_{134})}]$$
Replacing the $z_{ij}$ by $z^*_{ij}$ in the above expressions, we get the corresponding expressions for $[z^*_{13}z^*_{23}]$,$[z^*_{24}z^*_{14}]$,
$[z^*_{31}z^*_{41}]$ and $[z^*_{32}z^*_{42}]$.
By Proposition \ref{*prop} and the property of the triple ratio under the permutation of the indices, we obtain
$$
 z_{31}(1+z_{123})=z_{31}(1+z_{312})=z^*_{13}(1+\frac{1}{z_{134}})=z^*_{13}(1+z^*_{134}),
$$
and
$$
z_{32}\frac{z_{123}}{(1+z_{123})}=z_{32}\frac{1}{1+\frac{1}{z_{123}}}=z_{32}\frac{1}{1+\frac{1}{z_{231}}}=z^*_{23}\frac{1}{1+z_{324}}=
z^*_{23}\frac{1}{1+\frac{1}{z^*_{324}}}=z^*_{23}\frac{z^*_{243}}{(1+z^*_{243})}.
$$

Similarly, we have
$$ z_{42}(1+z_{142})=z^*_{24}(1+z^*_{243}),\; z_{13}(1+z_{134})=z^*_{31}(1+z^*_{123}),\; z_{24}(1+z_{243})=z^*_{42}(1+z^*_{142});
$$
and

$$
 z_{41}\frac{z_{142}}{(1+z_{142})}=z^*_{14}\frac{z^*_{134}}{(1+z^*_{134})},\; z_{23}\frac{z_{243}}{(1+z_{243})}=z^*_{32}\frac{z^*_{123}}{(1+z^*_{123})}, \;
 z_{14}\frac{z_{134}}{(1+z_{134})}=z^*_{41}\frac{z^*_{142}}{(1+z^*_{142})}.
$$

Now substituting the expressions above into the formula for the difference $\beta(T)-\beta(T^*)$, since the last 2 terms cancel out, we obtain
 $$
  \beta(T)-\beta(T^*)=2([-z_{123}]+[-z_{243}]+[-z_{134}]+[-z_{142}])-(\beta(T)-\beta(T^*)),
  $$
where we use the fact $[\frac{1}{t}]=-[t]$ in ${\mathcal{P}}(\bC)$, $\forall \, z\in \bC\setminus\{0,1\}$. Therefore,
$$
2(\beta(T)-\beta(T^*))=2([-z_{123}]+[-z_{243}]+[-z_{134}]+[-z_{142}]).
$$
Since ${\mathcal{P}}(\bC)$ is uniquely divisible,
$$ \beta(T)-\beta(T^*)=[-z_{123}]+[-z_{243}]+[-z_{134}]+[-z_{142}].$$

Case (II). $z_{12}z_{21}=1$ or $z_{34}z_{43}=1$. If both $z_{12}z_{21}$ and $z_{34}z_{43}$ are equal to $1$, then by (\ref{inv}), $$[z_{12}]+[z_{21}]=0=[z_{34}]+[z_{43}].$$
Hence $\beta(T)=0$. By Corollary \ref{cor:ij=kl}, $\beta(T^*)=0$. By (\ref{faceedge}) and (\ref{Cros}), we calculate that $$[-z_{123}]+[-z_{243}]+[-z_{134}]+[-z_{142}]=0.$$ Hence the proposition holds in this case.

Now suppose that $z_{12}z_{21}=1$ and $z_{34}z_{43}\ne 1$. Then $$\beta(T)-\beta(T^*)=[z_{34}]+[z_{43}]-[z^*_{12}]-[z^*_{21}].$$ Using the corresponding formulae in the case (I), we see that the proposition holds. We omit the details here since it is analogous to the case (I). For the case $z_{12}z_{21}\ne 1$ and $z_{34}z_{43}=1$, the proof is similar. The proposition is proved.
\EPf

\subsection{Remark}
By the Proposition, in $\mathcal{P}(\bC)$, $\beta(T)$ and $\beta(T^*)$  differ exactly by the face terms , so their volumes  differ by the sum of the values of the dilogarithm function $D$ at the negative face coordinate. In particular, their volumes are generically different.  This can be easily checked numerically.
%\end{remark}

\subsection{Remark}
A slightly different proof of the main proposition above is obtained by defining another set of coordinates
$w_{ij}=z_{ij}z_{ji}$, $1\leq i\ne j\leq 4$.  We clearly have $w_{ij}=w_{ji}$.  The map defined in this way is birational.  We have indeed
$$
z_{12} = w_{12}\frac{w_{13}w_{23}-w_{23}+1}{w_{12}w_{13}w_{23}+1}\ \ \ \ \  z_{21} = \frac{w_{12}w_{13}w_{23}+1}{w_{13}w_{23}-w_{23}+1}
$$
$$
 z_{34} = w_{34}\frac{w_{13}w_{14}-w_{14}+1}{w_{13}w_{14}w_{34}+1}\ \ \ \ \  z_{43} = \frac{w_{13}w_{14}w_{34}+1}{w_{13}w_{14}-w_{14}+1}.
 $$
The $w$-coordinates have the advantage that under the
duality $\Theta_4$ they have a very simple expression, namely, for $\{i,j,k,l\}=\{1,2,3,4\}$
$$
w^*_{ij}=w_{kl}.
$$
Observe also that the triple ratio has a simple formula in those coordinates
$$
-z_{21}z_{31}z_{41}= \frac{1}{w_{34}w_{41}w_{13}},
$$
and similarly for the other triple ratios.

To prove the main proposition we express the invariant $\beta(T)$ in those coordinates.  Using the proof of Proposition \ref{proposition:main}, we obtain
$$
2\beta(T)=
$$
$$
[w_{12}]+ [w_{21}]+[w_{34}]+ [w_{43}]-([w_{12}w_{13}w_{23}]+[w_{13}w_{34}w_{41}]+[w_{14}w_{42}w_{21}]+[w_{24}w_{43}w_{32}]) +R,
$$

where $R$ is the sum of the following terms (the identities below are, in fact, in $\mathcal{P}(\bC)$)
$$[z_{42}(1-z_{13}z_{23}z_{43})]+[z_{41}\frac{z_{13}z_{23}z_{43}}{z_{13}z_{23}z_{43}-1}]=[\frac{w_{12}w_{14}-w_{12}+1}{w_{12}w_{14}}]+[\frac{1}{w_{24}w_{12}-w_{24}+1}]$$
$$[z_{31}(1-z_{14}z_{24}z_{34})]+[z_{32}\frac{z_{14}z_{24}z_{34}}{z_{14}z_{24}z_{34}-1}]=[\frac{w_{12}w_{23}-w_{12}+1}{w_{23}w_{12}}]+[\frac{1}{w_{31}w_{12}-w_{31}+1}]$$
$$[z_{24}(1-z_{21}z_{31}z_{41})]+[z_{23}\frac{z_{21}z_{31}z_{41}}{z_{21}z_{31}z_{41}-1}]=[\frac{w_{34}w_{23}-w_{34}+1}{w_{23}w_{34}}]+[\frac{1}{w_{24}w_{34}-w_{24}+1}]$$
$$[z_{13}(1-z_{12}z_{32}z_{42})]+[z_{14}\frac{z_{12}z_{32}z_{42}}{z_{12}z_{32}z_{42}-1}]=[\frac{w_{34}w_{14}-w_{34}+1}{w_{34}w_{14}}]+[\frac{1}{w_{13}w_{34}-w_{13}+1}]$$

It is clear then that there is a symmetry under the duality.
%\end{remark}

\subsection{Examples}
Now we will give two geometric structures which can be naturally interpreted by the very generic tetrahedra of flags. One is from the the hyperbolic space $\mathbb H^3$, the other is from the $3$-dimensional spherical CR geometry.  The case of
$3$-dimensional flag manifolds structures is just a real version of the complex flags case.
We will also see under the duality $\Theta_4$, they have very simple formulae.
\begin{example}
\textbf{ Hyperbolic ideal tetrahedra.} Recall that an ideal hyperbolic tetrahedron is given by $4$ points on the boundary of $\mathbb H^3$,
i.e. a projective line $\mathbb{CP}^{1}$. It is known that the group of orientation-preseving isometries of $\mathbb H^3$ is $\PSL(2,\bC)$ and up to the action of this group, these points are in homogeneous coordinates $[0,1]$, $[1,0]$, $[1,1]$ and $[1,z]$ where $z\in \bC \setminus \{0,1\}$ is the cross-ratio of these four points.

Identifying $\bC^3$ with the Lie algebra $\mathfrak{sl}(2,\bC)$, we have the adjoint action
of $\PSL(2,\bC)$ on $\bC^3$ preserving the quadratic form defined by the determinant, which is
given in canonical coordinates by $xz-y^2$. The group $\PSL(2,\bC)$ preserves the
isotropic cone of this form. The projectivization of this cone is identified to $\mathbb{CP}^{1}$
via the Veronese map (in canonical coordinates):
\begin{eqnarray*}
 h_1 \: : \:  \mathbb{CP}^{1} & \to & \mathbb{CP}^{2} \\
 \left[x,y\right] & \mapsto &  [x^2,xy,y^2]
\end{eqnarray*}
The first jet of that map gives a map $h$ from $\mathbb{CP}^{1}$ to the space of flags $\Fl(V)$.
A convenient description of that map is obtained thanks to the
identification between $\bC^3$ and its dual given by the quadratic form. Denote $\langle ,\rangle$
the bilinear form associated to the determinant. Then we have
\begin{eqnarray*}
 h \: : \:\mathbb{CP}^{1} & \to & \Fl(V)\\
 p & \mapsto &  \left(h_1(p), \langle h_1(p), \cdot \rangle \right).
\end{eqnarray*}

Let $T$ be the tetrahedron $h([0,1])$, $h([1,0])$, $h([1,1])$ and $h([1,z])$. We call it a hyperbolic tetrahedron. By direct computations, we have the edge and face coordinates :
$$ z_{12}(T) = z_{21}(T)=z_{34}(T)=z_{43}(T)=z;$$
$$ z_{ijk}=1, \; \forall i,j,k \in \{1,2,3,4\}.$$
Hence, by Proposition \ref{4flags}, $T$ is very generic. By Proposition \ref{*prop}, we find $$z_{ij}(T^*)=z_{ij}, \; 1\leq i\ne j\leq 4.$$ That is, a hyperbolic tetrahedron is self-dual. Therefore, $\beta(T)=\beta(T^*)$. Conversely, given parameters satisfying the equations above, they define a unique hyperbolic tetrahedron.
\end{example}

\begin{example}
\textbf{The Spherical CR case.} Spherical CR geometry is modeled on the sphere ${\mathbb S}^3$ equipped with a natural $\mathrm{PU}(2,1)$ action.  More precisely, consider the group $\mathrm{U}(2,1)$ preserving the Hermitian form
$\langle z,w \rangle = w^*Jz$ defined on $V={\bC}^{3}$ by
the matrix
$$
J=\left ( \begin{array}{ccc}

                        0      &  0    &       1 \\

                        0       &  1    &       0\\

                        1       &  0    &       0

                \end{array} \right )
$$
and put:
$$
        V_0 = \left\{ z\in {\bC}^{3}-\{0\}\ \ :\ \
 \langle z,z\rangle = 0 \ \right\},
$$
$$
        V_-   = \left\{ z\in {\bC}^{3}\ \ :\ \ \langle z,z\rangle < 0
\ \right\}.
$$

Let $\pi :{\bC}^{3}\setminus\{ 0\} \rightarrow \P(V)=\mathbb{CP}^{2}$ be the
canonical projection.  We define
${\mathbb H}_{\bC}^{2} = \pi(V_-)$  the complex hyperbolic space and its boundary $\partial{\mathbb H}_{\bC}^{2}$ can be identified with ${\mathbb S}^3$, that is,
$$
\partial{\mathbb H}_{\bC}^{2} ={\mathbb S}^3= \pi(V_0)=\{ [x,y,z]\in \mathbb{CP}^{2}\ |\ x\bar z+ |y|^2+z\bar x=0\ \}.
$$

This defines a natural inclusion $h_1 \: : \:  {\mathbb S}^3 \mapsto  \mathbb{CP}^{2}$. The group of biholomorphic transformations of ${\mathbb H}_{\bC}^{2}$ is $\mathrm{PU}(2,1)$, the projectivization of $\mathrm{U}(2,1)$.  It acts on ${\mathbb S}^3$ by CR transformations. We
define $\bC$-circles as boundaries of complex lines in ${\mathbb H}_{\bC}^{2}$. Analogously,
$\bR$-circles are boundaries of totally real totally geodesic two dimensional submanifolds in ${\mathbb H}_{\bC}^{2}$.

Now an element $x\in {\mathbb S}^3$ gives rise to a flag  $(h_1(x),l_x)$, where $l_x$ is the unique complex line tangent to ${\mathbb S}^3$ at $x$. Therefore,we have a map
\begin{eqnarray*}
 h \: : \: {\mathbb S}^3 & \to & \Fl(V)\\
 x & \mapsto &  \left(h_1(x),l_x \right).
\end{eqnarray*}
Recall that $4$ distinct points in ${\mathbb S}^3$ are generic if any three of them are not contained in a $\bC$ circle. Given $4$ generic points $p_1$, $p_2$, $p_3$ and $p_4$, we have a tetrahedron of flags $T$ given by $h(p_1)$, $h(p_2)$, $h(p_3)$ and $h(p_4)$. We call it \emph{a spherical CR tetrahedron}. Then by definition its edge coordinates $z_{ij}$ are exactly those defined in \cite{F2}. In addition to the relations (\ref{Cros}), they satisfy another three equations:
\begin{eqnarray}\label{eq:cr}
z_{ij}z_{ji}=\overline {z}_{kl}\overline {z}_{lk}
\end{eqnarray}
with not all of them being real. By \cite{F3} page 411, the face coordinates are given by
$$
z_{ijk}=-z_{il}z_{jl}z_{kl}=\exp{2i\mathbb{A}(p_i, p_j, p_k)},
$$
where $\mathbb{A}(p_i, p_j, p_k)$ is the Cartan invariant.

Hence $T$ is very generic if and only if $\mathbb{A}(p_i, p_j, p_k)\ne \pm \frac{\pi}{2}$. This is equivalent to saying that $p_i$, $p_j$, $p_k$ are not on a $\bC$-circle. By \cite{F2} Proposition 4.2, we see that $T$ is very generic if and only if $p_1$, $p_2$, $p_3$ and $p_4$ are generic points on ${\mathbb S}^3$. Thus, a spherical CR tetrahedron $T$ is very generic. For the dual coordinates, we have
$$
z^*_{ij}=\overline{z}_{ij}.
$$
I.e. the dual of a spherical CR tetrahedron corresponds to the complex conjugation.

Note the condition $z_{ij}z_{ji}=\overline {z}_{kl}\overline{z}_{lk} \in \bR$ for all equations describes other configurations of four flags, not necessarily CR configurations, but including
the ones contained in an $\bR$-circle which coincide with the real hyperbolic ones with real cross ratio.  That is,
$$ z_{12}(T) = z_{21}(T)=z_{34}(T)=z_{43}(T)=x\in \bR\setminus \{0,1\}.$$
\end{example}

\section{Decorations of triangulations, $\mathcal{P}(\bC)$-valued invariant, and duality}
\subsection{Ideal triangulations by flag tetrahedra.}

We use the definition of an ideal triangulation of an oriented 3-manifold as it is used in computations with Snappea. It is a finite union of $3$-simplices $K= \bigcup_{\nu=1}^{N} T_\nu$ with face identifications (which are simplicial maps).  The quotient manifold $|K|-| K^{(0)} |$
can be retracted to a compact manifold with boundary when the vertices $| K^{(0)} |$ of the tetrahedra are deleted. The boundary has
two types, one coming from 2-skeletons on the boundary of $K$ and the other one being the links of $|K^{(0)}|$.
Observe that, sometimes, $K$ is called a quasi-simplicial complex (see \cite{N} section 4). It differs from a simplicial complex by the fact that there might be identifications on the boundary of a simplex.

Given such a compact oriented $3$-manifold $M$ with boundary, we call a triangulation $K$ as above an {\it{ideal triangulation} }of $M$.  But we will also consider a general ideal triangulation $K$ such that $|K|-| K^{(0)} |$ is not necessarily a manifold.

\begin{definition}
A $\it generic\   decoration$ of an ideal triangulation $K= \bigcup_{\nu=1}^{N} T_\nu$ is a map from the $0$-skeletons of the triangulation to $\Fl(V)$) such that, for each (oriented) $3$-simplex $T_\nu$, the $4$ flags assigned to its $4$ vertices form a generic tetrahedron of flags. So a decoration associates for each $3$-simplex $T_\nu$ a set of coordinates $z_{ij}(T_{\nu})$ as defined in section \ref{ss:coord-para}.
We denote a decorated ideal triangulation $K$ by $(K,z)$, where $z=(z_{ij}(T_{\nu}))_{\nu=1}^{N}$ stands for the set of coordinates associated for all $3$-simplices.
\end{definition}
 
From now on we fix $K=\bigcup_{\nu=1}^{N} T_{\nu}$ a decorated oriented ideal triangulation together with an ordering of the vertices of each $3$-simplex of $K$. Let $z_{ij}(T_{\nu})$ be the corresponding $z$-coordinates. We define the following compatibility conditions which imply the existence of a well defined representation of the fundamental group of the manifold $M$ into $\PGL(3,\bC)$ (cf. \cite{BFG} and \cite{F2} for the CR case).

\subsection{Consistency relations}\label{sec-consistency}
Let $F$ be a face ($2$-dim simplex) of $|K|$. Let $T$ and $T'$ be two $3$-simplices of $K$ with a common face $F$. Suppose that the vertices of $T$ (resp. $T'$) are $i,j,k,l$ (resp.$i,j,k,l')$ and that the face $F$ is $(ijk)$ which is oriented as a boundary of $T$. Note when the face $F$ is attached to $T$ or $T'$, it only changes its orientation. Hence we have the following for the face coordinates:

\medskip

\emph{(Face equations)}
Let $T$ and $T'$ be two tetrahedra of $K$ with a common face $(ijk)$ (oriented as a boundary of $T$), then $z_{ijk}(T)z_{ikj}(T')=1$.

\medskip

For a fixed edge $e\in |K|$ let $T_{\nu_1}, \cdots , T_{\nu_{n_e}}$ be the $n_e$ $3$-simplices in $K$ which contain an edge which projects onto $e\in |K|$ (counted with multiplicity).  For each $3$-simplex in $K$ as above we consider its edge $ij$ corresponding to $e$. We have the following for the edge coordinates:

\medskip

\emph{(Edge equations)} $z_{ij}(T_{\nu_{1}}) \cdots z_{ij}(T_{\nu_{n_e}})=z_{ji}(T_{\nu_1}) \cdots z_{ji}(T_{\nu_{n_e}})=1$.

\begin{definition} A $\emph{generic parabolic decoration}$ of an ideal triangulation $K=\bigcup_{\nu=1}^{N} T_\nu$ is a generic decoration $(K,z)$ satisfying the the consistency relations, i.e. the above face equations and edge equations. A $\emph{very generic parabolic decoration}$ of an ideal triangulation $K$ is a parabolic decoration such that for each $3$-simplex $T_\nu$ of $K$, the $4$ flags associated to its $4$ vertices form a very generic tetrahedron of flags.
\end{definition}

\begin{definition}
Let $K=\bigcup_{\nu=1}^{N} T_\nu$ be an ideal triangulation. We define $\text{Sol}(K)$ to be the set of all generic parabolic decorations of $K$. Since the face and edge equations are polynomials in $z_{ij}(T_\nu)$ with integer coefficients, $\text{Sol}(K)$ is a quasi-projective variety defined over $\mathbb{Z}$. Clearly the set of all very generic parabolic decorations is a Zariski open subset of $\text{Sol}(K)$.
\end{definition}

Given a decoration $(K,z)$, if for each $3$-simplex $T_\nu$, we replace the corresponding tetrahedron of flags by its complex conjugate, then we get a new decoration. We call it the complex conjugate of $(K,z)$ and denote it by $(K,\bar{z})$. Note that the coordinates of $(K,\bar{z})$ are exactly the complex conjugate of those of $(K,z)$. Since the face and edge equations are polynomials in $z_{ij}(T_\nu)$ with integer coefficients, if $(K,z)$ is parabolic, so is $(K,\bar{z})$. Analogously, given a very generic parabolic decoration $(K,z)$, if for each $3$-simplex $T_\nu$, we replace the corresponding tetrahedron of flags by its dual, then we get a new decoration. We call it the dual of $(K,z)$ and denote it by $(K,z^*)$. The coordinates of $(K, z^*)$ are related to those of $(K,z)$ by Proposition \ref{*prop}. Next we show that the dual of a very generic parabolic decoration is also very generic parabolic.

\begin{proposition}
Let $(K,z)$ a very generic parabolic decoration. Then $(K,z^*)$ is a very generic parabolic decoration.
\end{proposition}

\Pf Since the dual of a very generic tetrahedron of flags is very generic, it suffices to show that $(K,z^*)$ satisfies the face and edge equations. For the face equations, since $(K,z)$ satisfies them, it follows from Proposition \ref{*prop} $(1)$ that $(K,z^*)$ also satisfies. Let $e\in |K|$
be an edge. Suppose $T_{1}, \cdots , T_{n}$ are the $n$ $3$-simplices in $K$ which contain an edge $(ij)$ projecting onto $e$. Hence the vertices of these $3$-simplices can be given as follows:
$$
   T_1=(ijk_1k_2),\;T_2=(ijk_2k_3),\;\cdots,\; T_{n-1}=(ijk_{n-1}k_n),\;T_n=(ijk_nk_1).
$$
By Proposition \ref{*prop} $(2)$, we have
$$
z^*_{ij}(T_{1})\cdots z^*_{ij}(T_n)=z_{ji}(T_1)\frac{1+z_{jik_{2}}(T_{1})}{1+\frac{1}{z_{ijk_{1}}(T_{1})}}z_{ji}(T_2)\frac{1+z_{jik_{3}}(T_{2})}{1+\frac{1}{z_{ijk_{2}}(T_{2})}}
\cdots z_{ji}(T_n)\frac{1+z_{jik_{1}}(T_{n})}{1+\frac{1}{z_{ijk_{n}}(T_{n})}}.
$$
Since $T_{l-1}$ and $T_{l}$ have a common face $(ijk_{l})$, $2\leq l\leq n$, $T_n$ and $T_1$ have a common face $(ijk_1)$, by face equations, we obtain
$$
\frac{1}{z_{ijk_{l}}(T_{l})}=z_{jik_{l}}(T_{l-1}),\; 2\leq l\leq n,\;\; \text{and}\; \; \frac{1}{z_{ijk_{1}}(T_{1})}=z_{jik_{1}}(T_{n}).
$$
Since $(K,z)$ satisfied the edge equations, we have
$$
z^*_{ij}(T_{1})\cdots z^*_{ij}(T_n)=z_{ji}(T_1)\cdots z_{ji}(T_n)=1.
$$
Similarly, we see that
$$
z^*_{ji}(T_{1})\cdots z^*_{ji}(T_n)=1.
$$
Therefore, $(K,z^*)$ is a very generic parabolic decoration.
\EPf

By the above Proposition, we see the duality induces a birational involution on $\text{Sol}(K)$ which commutes with the complex conjugation. Hence $\text{Sol}(K)$ has a subgroup of birational automorphisms isomorphic to $\mathbb{Z}/2\mathbb{Z} \times \mathbb{Z}/2\mathbb{Z}$.

\begin{definition}\label{def beta}
Let $(K,z)\in \text{Sol}(K)$ be a parabolic decoration. We define the invariant
$$\beta(K,z)=\sum_{\nu=1}^{N}([z_{12}(T_\nu)]+[z_{21}(T_\nu)]+[z_{34}(T_\nu)]+[z_{43}(T_\nu)])\in \mathcal{P}(\bC).
$$
We define its volume to be $$ \text{Vol}(K,z)=\frac{1}{4}D(\beta(K,z)),$$ where $D$ is the Bloch-Wigner dilogarithm function.
\end{definition}

Now we have a map $$\beta: \text{Sol}(K)\rightarrow \mathcal{P}(\bC),\;\; (K,z)\mapsto \beta(K,z).$$

Recall we have an involution $\tau$ on the pre-Bloch group $\mathcal{P}(\bC)$ induced by complex conjugation. Let $(K,z)\in \text{Sol}(K)$ and $(K,\bar{z})$ be its complex conjugation. It is clear that
$$\beta(K,\bar{z})=\tau(\beta(K,z)).$$
Since $D(\bar{z})$=-$D(z)$ for any $z\in \bC$, we have
$$
\text{Vol}(K,\bar{z})=-\text{Vol}(K,z),\; \forall (K,z)\in \text{Sol}(K).
$$

Now we consider $\beta(K,z)$ and $\beta(K,z^*)$ for a very generic parabolic decoration. We have the following theorem:

\begin{theorem}\label{theorem:main}
Let $M$ be a $3$-dimensional manifold with an ideal triangulation $K= \bigcup_{\nu=1}^N T_\nu$. Assume $M$ has no boundary coming from $2$-skeletons on the boundary of $K$. Let $(K,z)\in \text{Sol}(K)$ be a very generic parabolic decoration and $(K,z^*)$ be its dual. Then, in ${\mathcal{P}}(\bC)$,
$$
\beta(K,z)=\beta(K,z^*).
$$
Hence $(K,z)$ and $(K,z^*)$ have the same volume.
\end{theorem}
\Pf
For each $3$-simplex $T_\nu$ of $K$, we still use $T_\nu$ to denote the very generic tetrahedron of flags associated to its four vertices. Then from the definition,
$$\beta(K,z)=\sum_{\nu=1}^N\beta(T_\nu),\;\;\beta(K,z^*)=\sum_{\nu=1}^N\beta(T^*_\nu).$$
By Proposition \ref{proposition:main}, for each $\nu$, the difference $\beta(T^*_\nu)-\beta(T_\nu)$ is determined by the four face coordinates of $T_\nu$. Since $M$ has no boundary coming from $2$-skeletons on the boundary of $K$, the faces of all $3$-simplex $T_\nu$ are glued in pair. By face equations, the corresponding face coordinates are inverse to each other. Recall we have in ${\mathcal{P}}(\bC)$, $[z^{-1}]=-[z]$, $\forall z\in \bC\setminus \{0,1\}$. Therefore,
$$\beta(K,z^*)-\beta(K,z)=\sum_{\nu=1}^N (\beta(T^*_\nu)-\beta(T_\nu))=0.$$
From this, it is clear that $(K,z)$ and $(K,z^*)$ have the same volume.
\EPf

\subsection{Remark}
 If $M$ has boundary coming from $2$-skeletons on the boundary of $K$, then $
\beta(K,z)$ and $\beta(K,z^*)$ are not necessarily equal in $\mathcal{P}(\bC)$. By the proof of the Theorem, their difference is determined by the face coordinates of the corresponding $2$-skeletons on the boundary of $K$. In this paper, we always consider manifolds satisfying the assumption of the above theorem, for example, $3$-dimensional hyperbolic manifolds with finite volume.

\subsection{Holonomy of a decoration}
\label{sec-holonomy}
Once a parabolic decoration is fixed, we obtain immediately a holonomy map as described in \cite{BFG}.  This gives rise to a representation from the fundamental group $\pi_1(M)$ to $\PGL (3,\bC)$. In fact, one obtains more structure attached to that representation, namely a decorated representation where the data of flags fixed by the holonomy at the cusps.

Note that if a very parabolic decoration $(K,z)$ gives rise to a representation $\rho: \pi_1(M)\rightarrow \PGL (3,\bC)$ via holonomy, then its dual $(K,z^*)$ gives $\rho^*:\pi_1(M)\rightarrow \PGL (3,\bC)$ which is $\rho$ followed by the Cartan involution $\Theta$ on $\PGL (3,\bC)$, i.e., transpose inverse of a matrix.

We also obtain a representation of the fundamental group of the boundary link at each vertex into $\PGL (3,\bC)$. We say the representation is $\sl unipotent$ if the generators of each boundary link are simultaneously represented in a parabolic subgroup.  That is, if for each cusp all generators are up to a global conjugation and a scalar multiplication of the form
$$ \left(
\begin{matrix}
1 & \star & \star\\
0 & 1 & \star \\
0 & 0 & 1
\end{matrix} \right ) .$$

\begin{theorem}
Let $M$ be a $3$-dimensional manifold with an ideal triangulation $K= \bigcup_{\nu=1}^N T_\nu$. Assume $M$ has no boundary coming from $2$-skeletons on the boundary of $K$. Let $(K,z)\in \text{Sol}(K)$ be a very generic parabolic decoration and $(K,z^*)$ be its dual. Suppose the boundary holonomy of $(K,z)$ is unipotent. Then, both $\beta(K,z)$ and $\beta(K,z^*)$ are in ${\mathcal{B}}(\bC)$, and
$$
\beta(K,z)=\beta(K,z^*).
$$
Hence $(K,z)$ and $(K,z^*)$ have the same volume.
\end{theorem}

\Pf  By the previous theorem, we have $ \beta(K,z)=\beta(K,z^*)$. By Theorem 5.5.1 in \cite{BFG} (see also \cite{FW}), $\beta(K,z)\in {\mathcal{B}}(\bC)$. Hence $\beta(K,z^*)=\beta(K,z)\in \mathcal{B}(\bC)$.
\EPf

\section{Applications}

\subsection{Self-duality and spherical CR decoration}
If a configuration of four flags has all coordinates satisfying $z_{ij}=z_{ji}=z_{kl}=z_{lk}$ we obtain from
Proposition \ref{*prop} that $z^*_{ij}=z_{ji}=z_{ij}$ for all $i,j$.  We conclude that decorated triangulations
satisfying these conditions for all simplices
are self-dual in the sense that they are fixed by the duality.  Observe that $\PGL(2,\bC)$
is self-dual as for any $g\in \PGL(2,\bC)$, ${(g^{-1}})^T$ is conjugated to $g$.  This explains the fact that
hyperbolic decorated triangulations are self-dual. In other words, we can not see this duality for representations to $\PGL(2,\bC)$.

\begin{definition}
A generic decoration $(K,z)$ is called spherical CR if every simplex of $K$ is decorated with spherical CR flags which are described in detail in section 3.7 Example 2.
\end{definition}

A similar observation holds for spherical CR decorated triangulations :

\begin{proposition} Suppose $(K,z)$ is a CR decorated triangulation. Then
$$
(K,z^*)=(K,\bar z).
$$
\end{proposition}
\Pf By (\ref{eq:cr}), a spherical CR decorated triangulation satisfies
$$
z_{ij}z_{ji}=\overline {z_{kl}z_{lk}}.
$$
By Corollary \ref{cor:ij=kl}, we have
$$
 z^*_{ij}z^*_{ji} = z_{kl} z_{lk}.
$$
The result follows.
\EPf

As a consequence of the proposition above and Theorem  \ref{theorem:main}, we obtain the following

\begin{theorem}\label{cr case}
Let $M$ be a $3$-dimensional manifold with an ideal triangulation $K= \bigcup_{\nu=1}^N T_\nu$. Assume $M$ has no boundary coming from $2$-skeletons on the boundary of $K$. Suppose the simplices are decorated with spherical
CR flags satisfying face pairing compatibilities. Then, in ${\mathcal{P}}(\bC)$,
$$
\beta(K,z)=\beta(K,\bar z).
$$
That is, $\beta(K,z)\in \mathcal{P}(\bC)^{+}$. Hence, the volume of $(K,z)$ is zero.
\end{theorem}

\Pf
By the Theorem \ref{theorem:main} and the previous proposition, we see $$\beta(K,z)=\beta(K,z^*)=\beta(K,\bar z).$$ Then $\beta(K,z)\in \mathcal{P}(\bC)^{+}$ follows from the definition of $\mathcal{P}(\bC)^{+}$ in section 2.1. By (\ref{dplus}), $$D(u)=0, \, \forall \, u\in \mathcal{P}(\bC)^{+},$$ the volume vanishes.
\EPf

\subsection{Remark}
Cusped Spherical CR manifolds satisfy $D(\beta(K,z))=0$ (this follows from
the fact that $H^3(\PU(2,1),\bR)=0$, see \cite{FW}). The above theorem is clearly a refined result.

\subsection{Cusped hyperbolic manifolds of finite volume}
The following theorem follows from an argument  in \cite{N} section 4.  See also \cite{BFG}.  It holds in more generality but, for simplicity, we will consider only a version for cusped hyperbolic manifolds.
\begin{theorem}
Let  $\rho : \pi_1(M)\rightarrow \PGL(3,\bC)$ be a decorated
representation of the fundamental group of a cusped hyperbolic manifold of finite volume. By this we mean an assignment of a flag to each cusp
such that each boundary holonomy preserves that flag.  Then,
there exists an ideal triangulation of $M$ and a parabolic decoration with holonomy precisely the given representation.
\end{theorem}

\begin{definition}\label{rho}
Let  $\rho : \pi_1(M)\rightarrow \PGL(3,\bC)$ be a decorated
representation of the fundamental group of a cusped hyperbolic manifold $M$ of finite volume. By the theorem above, we obtain a decorated quasi-simplicial complex $(K,z)$ for $\rho$. We define an invariant $\beta(\rho)$ for $\rho$ by
$$\beta(\rho)=\beta(K,z)\in \mathcal{P}(\bC).$$
\end{definition}

Note that $\beta(\rho)$ is well-defined since $\beta(K,z)\in \mathcal{P}(\bC)$ will not depend on the particular choice of triangulation because each triangulation is connected by 2-3 moves and the invariant is preserved by them as long as it is defined.  The subtle point is that if the triangulation is not sufficiently fine there will be configurations of flags which are not generic so that the coordinates $z$ are not defined.  A more efficient definition can be done (see \cite{BBI}) by extending the definitions
to degenerate configurations.

Now from Theorem 5.3 and Theorem \ref{theorem:main} we obtain our main application.
\begin{theorem}
Let $\rho : \pi_1(M)\rightarrow \PGL(3,\bC)$ be a decorated
representation of the fundamental group of a cusped hyperbolic manifold $M$ of finite volume. Let  $\rho^* : \pi_1(M)\rightarrow \PGL(3,\bC)$ be the dual decorated representation (transpose inverse).  Then $\beta(\rho)=\beta(\rho^*)$.
\end{theorem}

\Pf
Since $\beta(\rho)=\beta(K,z)$ and $\beta(\rho^*)=\beta(K,z^*)$, where $(K,z)$ is a very generic parabolic decoration for the decorated
representation $\rho$. Now it follows from Theorem \ref{theorem:main}.
\EPf


\begin{thebibliography}{ZZ99}


\bibitem[BFG]{BFG}N. Bergeron, E. Falbel, A. Guilloux ;
Tetrahedra of flags, volume and homology of SL(3). Geom. Topol. 18 (2014), no. 4, 1911-1971.

\bibitem[B1]{B1} S. Bloch ; Higher regulators, algebraic $K$-theory, and zeta functions of elliptic curves. CRM
    Monograph Series, 11. American Mathematical Society, Providence, RI, 2000.

\bibitem[BBI]{BBI} M. Bucher, M. Burger, A. Iozzi; Rigidity of representations of hyperbolic lattices $\Gamma < PSL(2,C)$ into $PSL(n, C)$. arXiv:1412.3428 [math.GT].

\bibitem[CURVE]{Cu} http://curve.unhyperbolic.org/

\bibitem[DGG]{DGG} T. Dimofte, M. Gabella, A. B. Goncharov ; K-Decompositions and 3d Gauge Theories. arXiv:1301.0192.

\bibitem[DS]{DS}  J. Dupont, C. H. Sah ; Scissors congruences. II. J. Pure Appl. Algebra 25 (1982), no. 2, 159--195.

\bibitem[F1]{F1} E. Falbel ; {Geometric structures associated to triangulations as fixed point sets of involutions}.
Topology and its Applications  154  (2007), no. 6, 1041-1052.

\bibitem[F2]{F2} E. Falbel ; { A spherical CR structure on the
complement of the figure eight knot with discrete holonomy}.
Journal of Differential Geometry 79 (2008) 69-110.

\bibitem[F3]{F3} E. Falbel ;
A volume function for Spherical CR Tetrahedra. Quarterly Journal of Mathematics 62  (2011),  no. 2, 397-415.

\bibitem[FGKRT]{FGKRT} E. Falbel, A. Guilloux, P.-V. Koseleff, F. Rouillier, M. Thistlethwaite ; Character varieties for $\SL(3,\bC)$: the figure eight knot.  Experimental Mathematics, 25 (2016), no. 2, 219–235.

\bibitem[FKR]{FKR} E. Falbel, P.-V. Koseleff, F. Rouillier ; Representations of fundamental groups of 3-manifolds into PGL(3;C): Exact computations in low complexity.  Geometriae Dedicata  177 (2015), 229–255.

\bibitem[FW]{FW} E. Falbel, Q. Wang ; A combinatorial invariant for Spherical CR structures. Asian J. Math. 17 (2013), no. 3, 391-422.

\bibitem[GGZ]{GGZ} S. Garoufalidis, M. Goerner, C. Zickert ; Gluing equations for PGL(n,C)-representations of 3-manifolds. Algebr. Geom. Topol.  15  (2015),  no. 1, 565–622.

\bibitem[GTZ]{GTZ} S. Garoufalidis, D. Thurston, C. Zickert ; The complex volume of SL(n,C)-representations of 3-manifolds. Duke Math. J.  164  (2015),  no. 11, 2099–2160.

\bibitem[Go]{Go} A. B. Goncharov ; Polylogarithms and motivic Galois groups. Motives (Seattle, WA, 1991), 43-96, Proc. Sympos. Pure Math., 55, Part 2, Amer. Math. Soc., Providence, RI, 1994.

\bibitem[G]{G}A. Guilloux ; Representations of 3-manifold groups in PGL(n,C) and their restriction to the boundary.  arXiv:1310.2907 [math.GT]

\bibitem[J]{J} H. ~Jacobowitz ; An Introduction to CR Structures.
Mathematical Surveys and Monographs {\bf 32}, American Math.\ Soc.\ (1990).

\bibitem[N]{N}W. D. Neumann ;  Extended {B}loch group and the {C}heeger-{C}hern-{S}imons class.  Geom. Topol. 8, (2004) 413–474.
	
\bibitem[NY]{NY} W. Neumann, Jun Yang ;
 Bloch invariants of hyperbolic 3-manifolds.
Duke Math. Journal 96 (1999) 25-59.

\bibitem[MFP]{MFP}  P. Menal-Ferrer, J. Porti ; Local coordinates for SL(n,C) character varieties of finite volume hyperbolic 3-manifolds, Ann. Blaise Pascal 19 no. 1 (2012), p. 107-122.

\bibitem[S1]{S1}A. A. Suslin ; $K_3$ of a field and the Bloch group.
Proc. of the Steklov Institute of Math. Issue 4 (1991), 217-238.

\bibitem[T]{T} W. Thurston ; The geometry and topology of 3-manifolds.
Lecture notes 1979.

\end{thebibliography}
\end{document}